\documentclass[twosided,reqno]{amsart}

\usepackage{amsmath}
\usepackage{amssymb}
\usepackage{amsthm}

\theoremstyle{theorem}
\newtheorem{theorem}{\scshape Theorem }[section]

\theoremstyle{definition}

\newtheorem*{remark*}{\scshape Remark}
\newcommand{\ma}{\mathbb}
\newcommand{\be}{\begin{equation}}
\newcommand{\ee}{\end{equation}}
\newcommand{\ben}{\begin{equation*}}
\newcommand{\een}{\end{equation*}}
\newcommand{\brn}{\begin{remark*}}
\newcommand{\ern}{\end{remark*}}
\newcommand{\fa}{\frac}
\newcommand{\la}{\label}
\newcommand{\U}{\sum_{n=0}^{\infty}}
\newcommand{\bt}{\begin{theorem}}
\newcommand{\et}{\end{theorem}}
\newcommand{\bi}{\binom}
\newcommand{\lp}{\left(}
\newcommand{\rp}{\right)}
\newcommand{\tn}{\frac{t^n}{n!}}
\newcommand{\ex}{\left(2\sum_{a=0}^{d-1}\fa{(-1)^a\chi(a)e^{at}}{e^{dt}+1}\right)^r e^{xt}}
\newcommand{\xq}{[x]_q}

\newcommand{\tq}{[2]_q^r}
\newcommand{\iq}{\lim_{q\to 1}}
\newcommand{\fq}{F_{q,\chi}^{(r)}(t,x)}
\newcommand{\fr}{F_{q,\chi}^{(r)}}
\newcommand{\eq}{E_{n,\chi,q}^{(r)}}
\newcommand{\en}{E_{n,\chi}^{(r)}}
\newcommand{\js}{j_1, \cdots ,j_r=0}
\newcommand{\is}{i_1, \cdots, i_r=0}
\newcommand{\ns}{n_1, \cdots, n_r=0}
\newcommand{\ms}{m_1, \cdots ,m_r=0}
\newcommand{\MS}{m_1+ \cdots +m_r}
\newcommand{\ty}{\infty}

\newcommand{\sr}{\sum_{l=1}^{r}}
\newcommand{\lr}{l_{q,r}}
\newcommand{\pl}{\Pi_{l=1}^{r}}

\numberwithin{equation}{section}

\begin{document}

\title{Some identities of symmetry for the generalized $q$-Euler polynomials}

\author{Dae San Kim}
\address{Department of Mathematics, Sogang University, Seoul 121-742, Republic of Korea.}
\email{dskim@sogang.ac.kr}

\author{Taekyun Kim}
\address{Department of Mathematics, Kwangwoon University, Seoul 139-701, Republic of Korea}
\email{tkkim@kw.ac.kr}

\maketitle

\begin{abstract}
By the symmetric properties of Drichlet's type multiple $q-l-$function, we establish various identities concerning the generalized higher-order $q$-Euler polynomials. Furthermore, we give some interesting relationship between the power sums and the generalized higher-order $q$-Euler polynomials.
\end{abstract}

\section{Introduction}

Let $\chi$ be a Dirichlet character with conductor $d\in \ma{N}$ with $d\equiv 1 \pmod{2}$. As is well known, the generalized higher-order Euler polynomials are defined by the generating function to be

\be\la{1}
\ex=\U\en(x)\tn.
\ee
When $x=0$, $\en=\en(0)$ are called the generalized Euler numbers attached to $\chi$ of order $r(\in\ma{N})$.\\
For $q\in\ma{C}$ with $|q|<1$, the $q$-number is defined by $\xq=\fa{1-q^x}{1-q}$.\\
Note that $\iq\xq=x$. In \cite{07}, Kim considered $q$-extension of generalized higher-order Euler polynomials attached to $\chi$ as follows:

\be\la{2}
\begin{split}
\fq&=\tq\sum_{\ms}^{\ty}(-q)^{\MS}\lp\Pi_{i=1}^{r}\chi(m_i)\rp e^{[\MS+x]_q t}\\
&=\U\eq(x)\tn.
\end{split}
\ee

Note that
\ben
\iq\fq=\ex.
\een

For $s\in\ma{C}$ and $x\in \ma{R}$ with $x\neq 0,-1,-2,\cdots,$ Kim defined Dirichlet-type multiple $q-l-$funtion which is given by

\be\la{3}
\begin{split}
l_{q,r}(s,x|\chi)&=\tq\sum_{\ms}^{\ty}\fa{(-q)^{\MS}\lp\Pi_{i=1}^{r}\chi(m_i)\rp}{[\MS+x]_{q}^{s}}\\
&=\fa{1}{\Gamma(s)}\int_{0}^{\ty}\fr(-t,x)t^{s-1}dt,~\textrm{(see \cite{07})}.
\end{split}
\ee

Applying the Laurent series and Cauchy residue theorem in (\ref{2}) and (\ref{3}), we get

\be\la{4}
\lr(-n,x|\chi)=\eq(x),~\textrm{where}~n\in\ma{Z}_{\geq 0}.
\ee

When $x=0,~\eq=\eq(0)$ are called the generalized $q$-Euler numbers attached to $\chi$ of order $r$. From (\ref{2}), we note that

\be\la{5}
\begin{split}
\eq(x)&=\sum_{l=0}^{n} \bi{n}{l}q^{lx}E_{l,\chi,q}^{(r)}\xq^{n-l}\\
&=\lp q^x E_{\chi,q}^{(r)}+\xq\rp^n,
\end{split}
\ee
with the usual convention about replacing $\lp E_{\chi,q}^{(r)}\rp^n$ by $\eq$, (see [1-13]).\\
In this paper, we investigate properties of symmetry in two variables related to multiple $q-l-$function which interpolates generalized higher-order $q$-Euler polynomials attached to $\chi$ at negative integers. From our investigation, we derive identities of symmetry in two variables related to generalized higher-order $q$-Euler polynomials attached to $\chi$. Recently, several authors have studied $q$-extensions of Euler polynomials due to T. Kim (see [1, 2, 3, 9, 10, 11, 12, 13]).

\section{Symmetry of $q$-power sum and the generalized $q$-Euler polynomials}

For $a,b \in \ma{N}$ with $a\equiv 1\pmod{2}$ and $b\equiv 1\pmod{2}$, we observe that

\be\la{6}
\begin{split}
\fa{1}{[2]_{q^a}^r}l_{q^a,r}&\lp s, bx+\fa{b}{a}(j_1+\cdots+j_r)|\chi\rp\\
&=[a]_q^s\sum_{\ns}^{\ty}\sum_{\is}^{db-1}\fa{(-1)^{\sr(i_l+n_l)}q^{a\sr(i_l+bdn_l)}\lp \Pi_{l=1}^{r}\chi(i_l)\rp}{\left[ab\lp x+d\sr n_l\rp +b\sr j_l+ a\sr i_l \right]_q^s}
\end{split}
\ee

From (\ref{6}), we have
\be\la{7}
\begin{split}
\fa{[b]_q^s}{[2]_{q^a}^r}&\sum_{\js}^{da-1}(-1)^{\sr j_l}\lp \pl \chi(j_l)\rp q^{b\sr j_l}l_{q^a,r}(s,bx+\fa{b}{a}\sr j_l|\chi)\\
&=[a]_q^s[b]_q^s \sum_{\is}^{db-1}\sum_{\js}^{da-1}\sum_{\ns}^{\ty}
\fa{(-1)^{\sr(i_l+n_l+j_l)}\lp \Pi_{l=1}^{r}\chi(j_l)\rp\lp \pl \chi(i_l)\rp}{\left[ab\lp x+d\sr n_l\rp +\sr (bj_l+ a i_l) \right]_q^s}\\
&\hspace{7cm}\times q^{\sr(ai_l+bj_l+abdn_l)}.
\end{split}
\ee

By the same method as (\ref{7}), we get

\be\la{8}
\begin{split}
\fa{[a]_q^s}{[2]_{q^b}^r}&\sum_{\js}^{db-1}(-1)^{\sr j_l}\lp \pl \chi(j_l)\rp q^{a\sr j_l}l_{q^b,r}(s,ax+\fa{a}{b}\sr j_l|\chi)\\
&=[a]_q^s[b]_q^s \sum_{\js}^{db-1}\sum_{\is}^{da-1}\sum_{\ns}^{\ty}
\fa{(-1)^{\sr(i_l+n_l+j_l)}\lp \Pi_{l=1}^{r}\chi(j_l) \chi(i_l)\rp}{\left[ab\lp x+d\sr n_l\rp +\sr (aj_l+ b i_l) \right]_q^s}\\
&\hspace{7cm}\times q^{\sr(aj_l+bi_l+abdn_l)}.
\end{split}
\ee

Therefore, by (\ref{7}) and (\ref{8}), we obtain the following theorem.

\bt\la{t1}
For $a,b\in \ma{N}$ with $a\equiv 1\pmod{2}$ and $b\equiv 1\pmod{2}$, we have
\ben
\begin{split}
&[2]_{q^b}^r [b]_q^s\sum_{\js}^{da-1}(-1)^{\sr j_l}\lp \pl \chi(j_l)\rp q^{b\sr j_l}l_{q^a,r}(s,bx+\fa{b}{a}\sr j_l|\chi)\\
&=[2]_{q^a}^r [a]_q^s\sum_{\js}^{db-1}(-1)^{\sr j_l}\lp \pl \chi(j_l)\rp q^{a\sr j_l}l_{q^b,r}(s,ax+\fa{a}{b}\sr j_l|\chi).
\end{split}
\een
\et

From (\ref{4}) and Theorem \ref{t1}, we obtain the following theorem.

\bt\la{t2}
For $n\geq 0$ and $a,b\in \ma{N}$ with $a\equiv 1\pmod{2}$ and $b\equiv 1\pmod{2}$, we have
\ben
\begin{split}
&[2]_{q^b}^r [a]_q^n \sum_{\js}^{da-1}(-1)^{\sr j_l}\lp \pl \chi(j_l)\rp q^{b\sr j_l} E_{n,\chi, q^a}^{(r)}(bx+\fa{b}{a}\sr j_l)\\
&=[2]_{q^a}^r [b]_q^n \sum_{\js}^{db-1}(-1)^{\sr j_l}\lp \pl \chi(j_l)\rp q^{a\sr j_l} E_{n,\chi, q^b}^{(r)}(ax+\fa{a}{b}\sr j_l).\\
\end{split}
\een
\et

By (\ref{5}), we easily get

\be\la{9}
\begin{split}
\eq(x+y)&=(q^{x+y}E_{\chi,q}^{(r)}+[x+y]_q)^n\\
&=(q^{x+y}E_{\chi,q}^{(r)}+q^x[y]_q+[x]_q)^n\\
&=\sum_{i=0}^{n}\bi{n}{i}q^{ix}(q^yE_{\chi,q}^{(r)}+[y]_q)^i[x]_q^{n-i}\\
&=\sum_{i=0}^{n}\bi{n}{i}q^{xi}E_{i,\chi,q}^{(r)}(y)[x]_q^{n-i}.\\
\end{split}
\ee

From (\ref{9}), we note that

\be\la{10}
\begin{split}
&\sum_{\js}^{da-1}(-1)^{\sr j_l}q^{b\sr j_l}\lp \pl \chi(j_l)\rp E_{n,\chi,q^a}^{(r)}(bx+\fa{b}{a}\sr j_l )\\
&=\sum_{\js}^{da-1}(-1)^{\sr j_l}q^{b\sr j_l}\lp \pl \chi(j_l)\rp \\
&\hspace{2cm}\times \sum_{i=0}^{n}\bi{n}{i}q^{ib\sr j_l}E_{i,\chi,q^a}^{(r)}(bx)\left[\fa{b(j_1+\cdots +j_r)}{a}\right]_{q^a}^{n-i}\\
&=\sum_{\js}^{da-1}(-1)^{\sr j_l}q^{b\sr j_l}\lp \pl \chi(j_l)\rp \\
&\hspace{2cm}\times \sum_{i=0}^{n}\bi{n}{i}q^{(n-i)b\sr j_l}E_{n-i,\chi,q^a}^{(r)}(bx)\left[\fa{b}{a}\sr j_l\right]_{q^a}^{i}\\
&=\sum_{i=0}^{n}\bi{n}{i} \lp\fa{[b]_q}{[a]_q}\rp^i E_{n-i,\chi,q^a}^{(r)}(bx)\\
&\hspace{2cm}\times \sum_{\js}^{da-1}(-1)^{\sr j_l}q^{b\sr (n-i+1)j_l}[j_1+\cdots +j_r]_{q^b}^{i}\\
&=\sum_{i=0}^{n}\bi{n}{i} \lp\fa{[b]_q}{[a]_q}\rp^i E_{n-i,\chi,q^a}^{(r)}(bx)S_{n,i,q^b}^{(r)}(ad|\chi),\\
\end{split}
\ee
where
\be\la{11}
S_{n,i,q}^{(r)}(a|\chi)=\sum_{\js}^{a-1}(-1)^{j_1+\cdots +j_r}\lp\pl \chi(j_l)\rp q^{\sr j_l(n-i+1)}\left[\sr j_l\right]_{q}^{i}.
\ee

From (\ref{10}) and (\ref{11}), we can derive the following equation.

\be\la{12}
\begin{split}
&[2]_{q^b}^r [a]_q^n \sum_{\js}^{da-1}(-1)^{\sr j_l}q^{b\sr j_l}\lp \pl \chi(j_l)\rp E_{n,\chi, q^a}^{(r)}(bx+\fa{b}{a}\sr j_l)\\
&=[2]_{q^b}^r \sum_{i=0}^{n}\bi{n}{i}[a]_q^{n-i}[b]_q^i E_{n-i,\chi,q^a}^{(r)}(bx)S_{n,i,q^b}^{(r)}(ad|\chi).
\end{split}
\ee

By the same method as (\ref{12}), we get

\be\la{13}
\begin{split}
&[2]_{q^a}^r [b]_q^n \sum_{\js}^{db-1}(-1)^{\sr j_l}q^{a\sr j_l}\lp \pl \chi(j_l)\rp E_{n,\chi, q^b}^{(r)}(ax+\fa{a}{b}\sr j_l)\\
&=[2]_{q^a}^r \sum_{i=0}^{n}\bi{n}{i}[b]_q^{n-i}[a]_q^i E_{n-i,\chi,q^b}^{(r)}(ax)S_{n,i,q^a}^{(r)}(bd|\chi).
\end{split}
\ee

Therefore, by  (\ref{12}) and (\ref{13}), we obtain the following theorem.

\bt\la{t3}
For $n\geq 0$ and $a,b\in \ma{N}$ with $a\equiv 1\pmod{2}$ and $b\equiv 1\pmod{2}$, we have
\ben
\begin{split}
&[2]_{q^b}^r \sum_{i=0}^{n}\bi{n}{i}[a]_q^{n-i}[b]_q^i E_{n-i,\chi,q^a}^{(r)}(bx)S_{n,i,q^b}^{(r)}(ad|\chi)\\
&=[2]_{q^a}^r \sum_{i=0}^{n}\bi{n}{i}[b]_q^{n-i}[a]_q^i E_{n-i,\chi,q^b}^{(r)}(ax)S_{n,i,q^a}^{(r)}(bd|\chi).
\end{split}
\een
\et

\brn
It is not difficult to show that
\be\la{14}
\begin{split}
&e^{[x]_q u}\sum_{\ms}^{\ty}q^{\sr m_l}(-1)^{\sr m_l}\lp \pl \chi(m_l)\rp e^{[y+\sr m_l]_q q^x (u+v)}\\
&=e^{-[x]_q u}\sum_{\ms}^{\ty}q^{\sr m_l}(-1)^{\sr m_l}\lp \pl \chi(m_l)\rp e^{[x+y+\sr m_l]_q (u+v)}.
\end{split}
\ee
\ern

Thus, by (\ref{14}), we get

\be\la{15}
\begin{split}
\sum_{k=0}^{m}\bi{m}{k}&q^{kx}E_{n+k,\chi,q}^{(r)}(y) [x]_q^{m-k}\\
&=\sum_{k=0}^{n}\bi{n}{k}q^{-kx}E_{m+k,\chi,q}^{(r)}(x+y) [-x]_q^{n-k}.
\end{split}
\ee


\end{document}